\title{Graph coloring with no large monochromatic components}

\newcommand{\cmt}[1]{\ifhmode\newline\fi{\sf *** \ \ #1 \\}}

\documentclass[11pt]{article}
\usepackage{a4}
\usepackage{latexsym}
\usepackage{epsfig}
\usepackage{amsmath}

\newtheorem{theorem}{Theorem}[section]

\newtheorem{lemma}[theorem]{Lemma}

\newtheorem{prop}[theorem]{Proposition}
\newtheorem{corol}[theorem]{Corollary}

\newcommand{\heading}[1]{\vspace{1ex}\par\noindent{\bf #1}}

\newcommand\mccsym{\mbox{\it mcc\/}}
\newcommand{\mcc}[1]{\mccsym_{#1}}
\newcommand{\mcctwo}{\mccsym_2}
\newcommand\cone{\mbox{\it cone\/}}

\makeatletter
\newcommand{\ProofEndBox}{{\ifhmode\unskip\nobreak\hfil\penalty50 \else
          \leavevmode\fi\quad\vadjust{}\nobreak\hfill$\Box$
            \finalhyphendemerits=0 \par}}%
\makeatother

\newcommand{\proofend}{\ProofEndBox\smallskip}

\newcommand\eps{\varepsilon}
\def\:{\colon}

\long\def\onefigure#1#2{
\begin{figure*}[tbp]
\begin{center}
#1
\end{center}
\caption{#2}
\end{figure*}
}

\newcommand{\lipefig}[2]  
{\onefigure{\mbox{\psfig{file=#1.eps}}}{\label{f:#1} #2} }

\newcommand\GG{{\mathcal{G}}}
\newcommand\LL{{\mathcal{L}}}
\newcommand{\ignore}[1]{}

\author{
{\sc Nathan Linial\thanks{Supported by the Israel Science Foundation.}}\\
{\footnotesize School of Computer Science and Engineering}\\[-1.5mm]
{\footnotesize  Hebrew University}\\[-1.5mm]
{\footnotesize  Jerusalem, Israel}\and
{\sc Ji\v{r}\'{\i} Matou\v{s}ek}\\
   {\footnotesize Department of Applied Mathematics and}\\[-1.5mm]
   {\footnotesize Institute of Theoretical Computer Science (ITI)}\\[-1.5mm]
   {\footnotesize  Charles University, Malostransk\'{e} n\'{a}m. 25}\\[-1.5mm]
   {\footnotesize  118~00~~Praha~1, Czech Republic}
\and
{\sc Or Sheffet\footnote{This paper is partially based on this author's
undergraduate \emph{Amirim} honours program project.}}\\
{\footnotesize School of Computer Science and Engineering}\\[-1.5mm]
{\footnotesize  Hebrew University}\\[-1.5mm]
{\footnotesize  Jerusalem, Israel}
\and
{\sc G\'abor Tardos}\thanks{Supported by NSERC grant 611470 and the Hungarian
 Foundation for Scientific Research Grant (OTKA) Nos. T037846, T046234,
 AT048826 and NK62321.}\\
{\footnotesize School of Computing Science}\\[-1.5mm]
{\footnotesize Simon Fraser University, Burnaby, BC, Canada}\\[-1.5mm]
{\footnotesize and R\'enyi Institute, Budapest, Hungary}
}
\date{}

\begin{document}
\maketitle
\addtolength\textwidth{-13mm}
\vskip1.5cm

\begin{abstract}
For a graph $G$ and an integer $t$ we let
$\mcc{t}(G)$ be the smallest $m$ such that
there exists a coloring of the vertices of $G$ by $t$ colors
with no monochromatic connected subgraph having more than $m$ vertices.
Let $\cal F$ be any nontrivial minor-closed family
of graphs. We show that 
$\mcctwo(G) = O(n^{2/3})$ for any $n$-vertex graph $G \in \cal F$.
This bound is asymptotically optimal and it is attained for
planar graphs. More generally, for every such $\cal F$, and every fixed $t$ 
we show that $\mcc t(G)=O(n^{2/(t+1)})$. 
On the other hand we have examples
of graphs $G$ with no $K_{t+3}$ minor and with $\mcc t(G)=\Omega(n^{2/(2t-1)})$.

It is also interesting to consider graphs of bounded degrees.
Haxell, Szab\'{o}, and Tardos proved
$\mcctwo(G) \leq 20000$ for every graph $G$ of maximum degree $5$.
We show that there are $n$-vertex $7$-regular graphs $G$
with $\mcctwo(G)=\Omega(n)$, and more sharply, for every $\eps>0$
there exists $c_\eps>0$ and $n$-vertex graphs of maximum
degree $7$, average degree at most $6+\eps$ for all subgraphs,
and  with $\mcctwo(G)\ge c_\eps n$.
For $6$-regular graphs it is known only that
the maximum order of magnitude of $\mcctwo$ is between $\sqrt n$
and $n$.

We also offer a Ramsey-theoretic perspective of the
quantity $\mcc{t}(G)$.
 
\end{abstract}

\section{Introduction}

In the classical graph coloring problem we assign a color to
each vertex so that no two vertices of the same color are adjacent.
In other words, each \emph{monochromatic connected component}
must be a single vertex. In the problems that
we study here, this requirement is relaxed and we
only demand that monochromatic connected components should
have small cardinality. Concretely, for a graph $G$ and
an integer $t$ we define $\mcc t (G)$ as the smallest integer $m$
such that the vertices of $G$ can be $t$-colored so that
no monochromatic connected component has cardinality exceeding $m$.
In particular, $\mcc t (G) = 1$ iff $G$ can be properly $t$-colored.

Here are some technicalities before we survey earlier work on this subject.
When we consider a graph $G$, $n$ always denotes the number of vertices.
For a set $S\subseteq V(G)$, we let $G[S]$ denote the subgraph of $G$ induced
by $S$. We use the standard asymptotic language and
conveniently ignore integrality issues in our computations.
Some of the examples we consider are line graphs. We recall
that the line graph of a graph $H$ is defined as $L(H)=
(E(H),\{\{e,e'\}: e\cap e'\ne\emptyset\})$. Clearly,
$\mcc t(L(H))$ is the smallest $m$ such that there is a
$t$-coloring of the edges of $H$ so that every
monochromatic connected subgraph of $H$ has $\le m$ edges.

The earliest reference investigating the parameter $\mcc t$
we are aware of 
is Kleinberg, Motwani, Raghavan, and Venkatasubramanian \cite{KMR+97},
where the question was motivated by a problem in computer
science concerning dynamically evolving databases.
Among others, the authors 
prove that $\mcc 3$ is unbounded for planar graphs, and 
that there is a constant $\eps>0$ such that for all (sufficiently
large) $d$, there are $d$-regular graphs $G$ with
$\mcc {\eps\sqrt d}(G)=\Omega(n)$.

Apparently independently, the possibility of bounding
$\mcc{t}(G)$  by a constant for graphs of bounded degree
has been investigated by graph theorists. The results
concern mainly the case $t=2$.
It is easy to see that $\mcctwo(G)\le 2$ for any graph $G$
of maximum degree 3.
Alon, Ding, Oporowski, and Vertigan \cite{ADO+03}
proved that $\mcctwo(G) \leq 57$ for every graph $G$ of maximum degree 4.
Haxell, Szab\'{o}, and Tardos \cite{HST03}
improved this  to $\mcctwo(G) \leq 6$ and
proved that $\mcctwo(G) \leq 20000$ for every graph $G$ of maximum degree $5$.
On the other hand, Alon et al.\ \cite{ADO+03} constructed
6-regular graphs $G$ with  $\mcctwo(G)$ arbitrarily large.
For graphs $G$ of maximum degree 3 it was also shown in \cite{BS05}
that they admit two-coloring where one color induces an independent  set,
while the other color induces components of size at most 189.
Earlier work on this subject \cite{DOS+}, \cite{JW96} mainly focused
on more specific questions concerning line graphs of 3-regular graphs.
These investigations culminated in \cite{Tho99} showing that the
edges of every 3-regular graph can be 2-colored so that each
monochromatic component is a path of length at most 5.

We should also mention that there is a fairly rich literature that
deals with the notion of $t$-vertex coloring where each
monochromatic connected component has a small {\em diameter},
see e.g.~\cite{LS93}. This line of research originated
in the field of distributed computing. 

Here are the main results of the present paper:

\begin{theorem}\label{t:fixedminor}
For every planar graph $G$ we have $\mcctwo(G)=O(n^{2/3})$.
More generally, for every nontrivial minor-closed
family of graphs $\cal F$ there exists
a constant $C=C_{\cal F}$ such that if $G$ belongs to $\cal F$,
then $\mcctwo(G)\le C n^{2/3}$. This bound is
tight even for planar graphs.

For every fixed integer $t\ge 2$ and every $G \in \cal F$
as above, there holds $\mcc t (G)=O(n^{2/(t+1)})$.
On the other hand, for every $t$ there exist graphs with no $K_{t+3}$
minor and with $\mcc t (G)=\Omega(n^{2/(2t-1)})$.

If the tree-width of $G$ is bounded by a constant (equivalently,
if $G$ excludes a fixed \emph{planar} minor), then
$\mcc t(G)=O(n^{1/t})$. This bound is asymptotically
optimal for every fixed $t$.

\end{theorem}

\begin{theorem}\label{t:6pluseps}
For every $\eps>0$ there exists a constant $c_\eps>0$
and arbitrarily large graphs $G$ so that
\begin{itemize}
\item
Every vertex in $G$ has degree at most $7$,
\item
Every subgraph of $G$ has average degree at most $6+\eps$,
\item
$\mcctwo(G)\ge c_\eps n$.
\end{itemize}
\end{theorem}

As already mentioned, the questions we consider here have
independently originated in computer science and in graph theory.
Graph coloring is, of course, one of the most fascinating parts
of graph theory. Due to its great significance and the famous
open questions about it, many different variations on the basic theme
are being investigated (see, e.g.,~\cite{JT95}) and the present
problems can be viewed as part of this ongoing research effort.

What is less obvious is the connection between
the graph invariants we consider here
and Ramsey Theory. The usual perception
is that Ramsey-type theorems express the fact that
large systems necessarily contain ``highly regular islands". We
suggest that many Ramsey-type results can
be viewed as ``sum theorems". Specifically, let $\cal G$ be a class
of graphs closed under taking subgraphs. Given a graph $G$, we ask
for the smallest number of members in $\cal G$  whose union is $G$.
Thus when $\cal G$ consists of all graphs without a $k$-clique
we encounter the classical Ramsey problem. When it is the class of
all graphs not containing a given subgraph we recover the
so-called Graph Ramsey Theory. Finally if $\cal G$ contains all
graphs in which each connected component has small cardinality
we arrive at our present problem. We believe that this
perspective deserves further research.
Needless to say, this concept extends beyond graphs and unions. Other
mathematical objects and other appropriate operations can
be considered. Given an object $G$ and a class $\cal G$, one seeks
the most economical way of expressing $G$ as a
``sum'' of members in $\cal G$.

There are certain classes of graphs for which the study
of coloring with small monochromatic connected components
is particularly interesting.
Let $D^d_m$ denote the $d$-dimensional grid with all diagonals;
that is, $D_m^d$ is the graph with vertex set $\{1,2,\ldots,m\}^d$
where the two vertices $u,v$ are adjacent if
$\|u-v\|_\infty=\max_{i}|u_i-v_i|\le 1$.
The study of $\mccsym$ for this graph and its relatives leads
to very interesting problems which bring together combinatorics
geometry and topology. The $d$-dimensional
version of the well-known HEX lemma
(see Gale~\cite{G79} and Linial and Saks \cite{LS93})
implies that $\mcc d (D_m^d)\ge m$
(and a simple coloring shows an $O(m)$ upper bound for every
fixed $d$). More colors allow for constant size components, i.e., we have $\mcc
{d+1}(D_m^d)=O(1)$ for every fixed $d$.
For two colors, Matou\v{s}ek and P\v{r}\'\i v\v etiv\'y
\cite{MP} proved $\mcctwo(D_m^d)\ge m^{d-1}-d^2 m^{d-2}$,
which nearly matches the obvious upper bound (layer-by-layer
$2$-coloring) of $m^{d-1}$. The behavior of $\mcc t (D_m^d)$
for $3\le t<d$ is still unknown and remains an intriguing open problem.

\section{Excluded minors, separators, and coloring}

In this section we prove all the 
upper bounds in Theorem~\ref{t:fixedminor}.

A subset $C\subseteq V$ of the vertex set of a graph $G=(V,E)$
is called a \emph{separator} if no component of $G[V\setminus C]$
has more than $\frac 23|V|$ vertices (the choice of the constant
$\frac 23$ is somewhat arbitrary). By the well-known planar separator
theorem of Lipton and Tarjan \cite{lip:tar}, every planar
graph $G$ has a separator with at most $O(\sqrt n\,)$ vertices.
More generally, for every $h$-vertex graph $H$, 
every $G$ containing no minor isomorphic
to $H$ has a separator with at most $h^{3/2}\sqrt n$ vertices \cite{AST}.
We will also need that graphs of bounded tree width  have
bounded size separators. In particular,
a graph of  tree width $w$ has a separator
of size $w+1$ as stated, e.g., 
in Remark~2 of \cite{Brucereed}. (We haven't found
the proof of this very simple fact anywhere, but a constant size
separator follows from the combination 
of the analogous result on branch width
in Lemma~3.1 of \cite{GM13} 
and the connection between branch width and tree width
as stated in \cite{GM10}.)

In view of these results,
the upper bounds in Theorem~\ref{t:fixedminor}
all follow from 
the following proposition.

\begin{prop}\label{p:separator}
 Let $\GG$ be a class of graphs closed under taking
induced subgraphs such that every $G\in \GG$
has a separator with at most $Kn^\gamma$ vertices, where $K$
and $\gamma\in [0,1)$ are constants depending only on $\GG$.
Then for every $G\in \GG$ we have
$$
\mcctwo(G)= O(n^{1/(2-\gamma)}),
$$
and more generally,
$$
\mcc t (G) = O(n^{1/(t-(t-1)\gamma)}),
$$
where the hidden constant of proportionality depends on
$K$, $\gamma$, and (in the second case) on $t$.
\end{prop}

\heading{Proof. } First we deal with the special case $t=2$.
Given an $n$-vertex $G$, let $n_0:=\lfloor{n^{1/(2-\gamma)}}\rfloor$
be a threshold parameter. We present a simple algorithm
producing a $2$-coloring of $V(G)$. The algorithm maintains
a list $\LL$ of induced subgraphs of $G$, which is initialized
to $\LL := \{G\}$, and a set $S$ of vertices,
initialized to $\emptyset$. While $\LL$ contains at least one graph
with more than $n_0$ vertices, we select one such graph $G_i\in\LL$
arbitrarily, we remove it from $\LL$,
we find a separator $C_i$ of $G_i$ 
of size at most $K |V(G_i)|^\gamma$, we set $S:=S\cup C_i$,
and we add all of the components of $G_i[V(G_i]\setminus C_i]$ to $\LL$.
The algorithm ends when $\LL$ contains only graphs of size at most
$n_0$; at this moment, we color the vertices of $S$ blue
and all remaining vertices (i.e., the vertices of all graphs 
in $\LL$) red and we finish. The algorithm is illustrated
in Fig.~\ref{f:mcc2}.

\lipefig{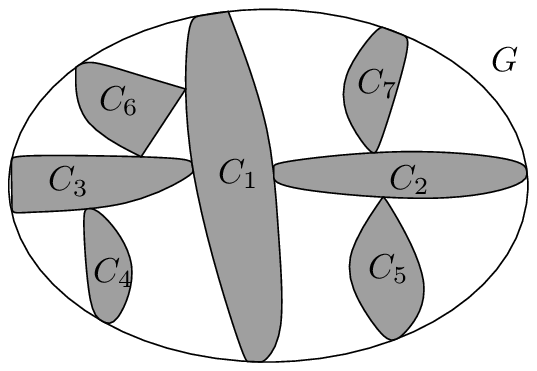}{The recursive algorithm for $2$-coloring via separators.}

By construction, no red component in this coloring has more than
$n_0$ vertices, and it suffices to show that $|S|=O(n_0)$
at the end of the algorithm 
($S$ can form a single blue component at worst).
We use the following
charging scheme: Whenever we color a separator $C_i$
in a graph $G_i$ blue, we let each vertex of $G_i$
pay $K|V(G_i)|^{\gamma-1}$ units. Since $|C_i|\le
K|V(G_i)|^\gamma$, the total paid by all vertices of $G_i$ at this step
is at least $|C_i|$. Now we consider an individual vertex
$v\in V(G)$ and we bound the total charge paid by it throughout the
whole algorithm. There may be several successive charges, since $v$ first
pays as a vertex of $G$, and then possibly as vertex of
some of the $G_i$. Let $x_j$ be the amount paid when $v$ is charged
the $j$th time, $j=1,2,\ldots,q$. 
We observe that
$x_q\le Kn_0^{\gamma-1}$ (since only graphs $G_i$
with at least $n_0$ vertices get partitioned), and that
$x_{j-1}\le (2/3)^{1-\gamma}x_{j}$ for all $j$, since the component
of $G_i$ containing $v$ always has at most two-third of the
vertices of $G_i$. Hence the $x_j$ are bounded
from above by a decreasing
geometric series, and thus the total charge paid by $v$
is $O(n_0^{\gamma-1})$. So $|S|\le O(nn_0^{\gamma-1})=O(n^{1/(2-\gamma)})$.

\medskip

Next, we consider the case of $t>2$ colors. We proceed
by induction on $t$, assuming that for every $n$-vertex
graph $G\in \GG$ we can construct a coloring with $t-1$ colors
witnessing $\mcc{t-1}(G) = O(n^{1/(t-1-(t-2)\gamma)})$.

For the induction step, we consider a $G\in \GG$ and we apply to it
the algorithm above with the following modifications:
This time we let the threshold be $n_0:=\lfloor n^{1/(t-(t-1)\gamma)}\rfloor$,
and at the end, we color $G[S]$ by $t-1$ colors
using the inductive assumption, while the vertices not belonging to $S$
get color $t$. 

\begin{sloppypar}
By the above analysis, we have $|S|=
O(nn_0^{\gamma-1})$, and by induction, the monochromatic components
in colors $1$ through $t-1$ have size at most
$O(|S|^{1/(t-1-(t-2)\gamma)})=O(n^{1/(t-(t-1)\gamma)})$. 
The components in color $t$ have size at most $n_0$ by construction.
This finishes the proof of Proposition~\ref{p:separator}.
\proofend
\end{sloppypar}

\section{Lower bounds for planar graphs and for excluded minors}

To prove that the bound in 
Theorem~\ref{t:fixedminor} for $\mcc2(G)$ is tight for planar
graphs $G$, 
we construct planar graphs $G$ with $\mcctwo (G)=\Omega(n^{2/3})$.
For every integer $k$ we construct $G=G_k$ on $n=2k^3+1$ vertices,
as indicated in Fig.~\ref{f:mcc3}.
\lipefig{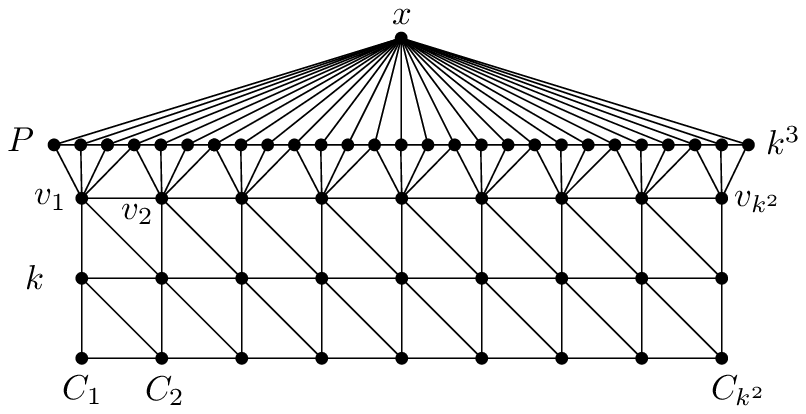}{The lower bound construction for planar graphs,
drawn for $k=3$.}
This  $G$ is constructed from a $k$ by $k^2$ grid $Z$, a path $P$
of $k^3$ vertices,
 and an extra vertex $x$. For $1\le i\le k^2$ we
denote the $i$th column of $Z$ by $C_i$, and we let $v_i$ be the
top vertex of $C_i$. 
We break the path $P$ into consecutive intervals
$I_1,\ldots,I_{k^2}$ of $k$ vertices each and we
connect the vertices of $I_i$
with $v_i$. We let $R_i=C_i\cup I_i$ and call this a 
\emph{rib} of $G$.
We connect $x$ with all vertices in $P$. Finally,
we add diagonals to all quadrilateral faces of
the planar graph constructed so far, so that 
it becomes a \emph{triangulated polygon}, that is,
a planar graph where all faces except possibly for one
are triangles.

Our main tool is the following lemma about triangulated
polygons (a very similar lemma
appears in \cite{MP}). 
For a set $S$ of vertices in a graph
we denote by $\partial S$ the set of vertices that are not in $S$ but
have a neighbor in $S$.

\begin{lemma}
\label{vertex_boundary}
Let $G$ be a triangulated polygon
and let $S\subseteq V(G)$ be such that
$G[S]$ is connected. 
Suppose that two vertices $u, v \in \partial S$ are not separated
by $S$ in $G$. Then there is a path between $u$ and $v$
that is entirely included in $\partial S$.
\end{lemma}

\heading{Proof.} This is a simple consequence of the planar HEX
lemma. Since $G[S]$ is connected and $u,v\in\partial S$,
 there is an $u$-$v$ path $P_1$ with all internal
vertices in $S$. Since $S$ doesn't separate $u$ from $v$,
there is another $u$-$v$ path $P_2$ that avoids $S$.
We consider subgraph of $H$ consisting of the
cycle $P_1\cup P_2$  plus the part of $G$ 
that triangulates the interior of this cycle
(assuming that the single non-triangular face 
of $G$ is the outer face). We add two new vertices $z$ and $t$
and we connect $z$ to all vertices of $P_1$ and 
$t$ to all vertices of $P_2$. The resulting
graph $H$ is a triangulation of the cycle $uzvt$.

We color blue all vertices in $\partial S$
including $u$ and $v$, and we 
color red all other vertices of $H$ including $z$ and $t$.
By the HEX lemma (as stated, e.g., in \cite{MN98})
we have either a blue $u$-$v$ path, which is what we want,
or a red $z$-$t$ path $Q$. We want to exclude the latter
possibility. 
Let us imagine that we follow the red path $Q$ 
from $t$ to $z$ and we watch the distance to $S$
in $H$. Since $t$ is adjacent only to $P_2$, whose vertices
are not in $S$, initially at $t$ this distance is at least $2$.
On the other hand, since the red vertices connected to $z$
are inner vertices of $P_1$ and thus in $S$, 
the penultimate vertex
in $Q$ is in $S$. Consequently, there is a vertex in
$Q$ at distance $1$ from $S$, but such a vertex was colored 
blue---a contradiction.
\proofend

\medskip

We need the following consequence of Lemma~\ref{vertex_boundary}:

\begin{corol}
Consider a red-blue vertex coloring of the graph $G=G_k$,
where $x$ is red.
Let $S$ be the connected component of 
$x$ in the red subgraph. If there is a connected
component of $G \setminus S$ containing at least
$r$ $S$-free ribs (i.e., ribs with no vertex in $S$),
then $G$ has a blue connected subgraph with
at least $rk$ vertices.
\end{corol}

\heading{Proof.}
If a rib $R_i$ is $S$-free, then the $k$ vertices in 
the interval $I_i$
are contained in $\partial S$ and are therefore blue.
For any two $S$-free ribs $R_i$, $R_j$ contained in
the same connected component of $G \setminus S$,
we choose vertices $u \in I_i$ and $v \in I_j$
(arbitrarily). The previous lemma now shows that
$u$ and $v$ are connected by a blue path.
\proofend

\medskip

We can now show that any two-coloring of $G$ has 
a monochromatic connected
subgraph of at least $k^2/2$ vertices.
As in the corollary, we assume $x$ red and we let
$S$ be the connected component of
$x$ in the red subgraph. We may assume $|S|\le k^2/2$,
for otherwise, we have a large red component.
Hence there are at least $k^2/2$ $S$-free ribs.
We want to show that at least $k/2$ of them
are in the same connected component of
$G \setminus S$; then we will be done by the corollary.

Since  $|S| < k^2$, at least one of the
$k$ \emph{rows} of the grid $Z$ contains fewer than
$k$ vertices from $S$. 
It follows that the $S$-free ribs live in at most $k$
connected components of $G \setminus S$. So there must be at least
$k/2$ of them in the same connected component 
as claimed.
\proofend
\bigskip

For $m\ge 1$ and 
a graph $G$ let $\cone(mG)$ be a graph constructed by
taking  $m$ isomorphic and pairwise disjoint copies 
of $G$ and connecting all of their vertices
to an additional new vertex, called the \emph{apex}.

\begin{lemma}
Let $t\ge 1$, let $G$ be a graph, and let $m=\mcc t(G)$.
 Then $\mcc{t+1}(\cone(mG))\ge m$.
\end{lemma}

\heading{Proof.} 
Let us consider a coloring of the vertices of 
$\cone(mG)$ with $t+1$ colors, and let us assume that
the apex has color $t+1$. Clearly, all
vertices of color $t+1$
form a connected subgraph, so if there are at least $m$
of them we have our large monochromatic connected subgraph. 
Otherwise, one of
the copies of $G$ lacks color $t+1$ and the claim follows.
\proofend
\medskip

Notice that as we pass from $G$ to $\cone(mG)$,
the number of vertices grows (approximately) $m$ times,
but other parameters grow slowly:
the tree width grows by at most
one, the size of the largest clique minor grows by
one, and if $G$ is outerplanar, then
$\cone(mG)$ is planar.

To prove the statement of Theorem~\ref{t:fixedminor}
about the existence of $K_{t+3}$ minor free graphs with 
high $\mcc t$, we simply take the planar
graph $G_k$ constructed at the beginning 
of this section and we apply the above lemma
$t-2$ times with $m=k^2/2$. 
Then $\mcc t$ of the resulting graph
is at least $k^2/2$, and the number of vertices
is $O(k^{2t-1})$.

To prove the similar statement 
about constant tree-width graphs,
we need a different base graph: 
let $F_k$ be the ``fan'' consisting of 
a $k$-vertex path and an 
additional vertex adjacent to all vertices of this
path. Clearly, 
$F_k$ is an outerplanar graph of tree width $2$.
A straightforward computation shows that 
$\mcc2(F_k)=\Theta(\sqrt k\,)$ (this also appears in
\cite{ADO+03}). 
Applying the above
lemma to $F_k$ $t-2$ times 
with $m=\mcc2(F_k)$, we obtain a graph of tree width
at most $t$ on $O(k^{t/2})$ vertices with
$\mcc t = \Omega(\sqrt k)$.
This finishes the proof of Theorem~\ref{t:fixedminor}.
\proofend

\medskip

The case $t=3$ in the just finished proof yields $n$-vertex
planar graphs with $\mcc 3$ at least $\Omega(n^{1/3})$.
 These very graphs were used
in \cite{KMR+97} to show that $\mcc3$ is not bounded for planar graphs.

\section{Edge expansion and degree \boldmath $6+\eps$}

In this section we prove Theorem~\ref{t:6pluseps}.  
The graphs we construct are line graphs $G=L(H)$.
So the property that we need is that in every $2$-coloring of $E(H)$
there are monochromatic connected components containing a positive fraction of the edges of $H$. 
To this end, it suffices to show that
small subgraphs of $H$ have small average degrees,
as the next observation shows:

\begin{lemma}\label{local_density}
Let $H$ be a graph with average degree $\overline d$. 
Suppose that every subgraph on $p$ or fewer vertices in $H$ 
has average degree strictly smaller than
$\overline d /t$. Then $\mcc t (L(H))\ge p$.
\end{lemma}

\heading{Proof. }
If $F \subseteq E(H)$ is the largest color class in a 
$t$-coloring of $E(H)$, 
then the graph $(V(H),F)$ has average degree at least 
$\overline d/t$. Consequently, some connected component of 
$(V(H),F)$ has average degree at least $\overline d/t$.
Such a connected component must 
have more than $p$ vertices, and thus at least $p$ edges.
\proofend

\medskip

In the proof of Theorem~\ref{t:6pluseps}, we will use
a suitable random graph (with bounded
vertex degrees) for $H$. The basic idea 
is that random graphs are typically very good \emph{expanders}.
This means that every sufficiently small set $S\subset V(H)$
has many neighbors, hence most of the edges incident to $S$
leave $S$, and consequently, the average degree
of $H[S]$ is small. More precisely, it turns out that
if $|S|$ is sufficiently small and if $H[S]$ is connected,
then this subgraph is nearly a tree in the sense that 
the average degree in $H[S]$ is just a little 
bigger than~$2$. Here sufficiently small means that
$|S| < \beta |V(H)|$, where $\beta>0$ is a suitable small constant
depending on how close we want to get to average degree $2$.

A result about this almost-tree behavior of small sets
in random \emph{regular} graphs
appears explicitly in \cite{HLW}. In particular,
Theorem~4.16, part (1) in \cite{HLW} tells us that
for every $d\ge 3$ and every $\delta>0$ there exists $\beta=
\beta(d,\delta)>0$
such that almost every $d$-regular graph $H$ on $m$ vertices
has average degree of $H[S]$ at most $2+\delta$
for all $S$ with at most $\beta m$ vertices. 
This, together with Lemma~\ref{local_density}, immediately
yields the following weaker analogue of Theorem~\ref{t:6pluseps}:
{\em There exist arbitrarily large
$8$-regular graphs $G$ with $\mcc 2(G)=\Omega(n)$.}
Indeed, we choose $H$ as a $5$-regular graph on $m$
vertices (thus, $G=L(H)$ is $8$-regular)
satisfying the conclusion of the statement
quoted above with $\delta=0.4$. Then for $\beta=\beta(5,0.4)$,
every $S\subseteq V(H)$ with at most $\beta m$ vertices
induces a subgraph average degree at most $2.4$, and 
hence Lemma~\ref{local_density} with $t=2$
and $\overline d=5$ shows $\mcc 2(L(H))\ge \beta m$.

In order to lower the maximum degree of $G$ to $7$ and
the maximum average degree to $6+\eps$, we will
use a random $H$ where a small fraction of vertices
have degree $5$, all others have degree $4$, and
no two degree-$5$ vertices are connected.

\heading{The random graph model. } 
It is easier to deal with random \emph{bipartite} graphs.
Most of the literature in this area deals with regular
random graphs, but we need a suitable mixture
of vertex degrees, and so we prescribe the degree
individually for each vertex.  That is, we 
have two disjoint sets $A$ and $B$ of vertices,
$|A| + |B|=m$,
and for every $v\in A\cup B$ we specify a number
$d(v)\in \{1,2,\ldots,D\}$, where $D$ is a constant.
These degrees and the sizes of $A$ and $B$ are related
by the condition $d(A)=d(B)$, where 
we use the notation $d(S)=\sum_{v\in S} d(v)$.
Hence $|A|=c_A m$, $|B|=c_B m$ for constants $c_A,c_B$.
To generate the random graph, every vertex $v$
starts with $d(v)$ ``half-edges'', and then 
the half-edges of all vertices in $A$ are matched at
random to the half-edges of all vertices in $B$
(this is a \emph{configuration model} of generating
random bipartite graphs). We note that the resulting $H$
may have multiple edges, but the line graph $L(H)$
we are interested in is still a simple graph.\footnote{It is
well known that for random regular graphs of fixed degree,
the configuration model yields a simple graph with
probability bounded away from $0$, and consequently,
any property that holds almost surely in the configuration
model also holds for almost all simple regular graphs
of the given degree; see, e.g., \cite{JLR}.
By slightly modifying the proof of this fact,
we could also get a similar
result  for our model with mixed degrees, and hence
have $H$ simple.}

The following lemma speaks about number of \emph{vertices}
adjacent to $S$; the number of edges is then obtained
as a simple consequence.

\begin{lemma}\label{sparse}
Let $D$ be a fixed integer. Then 
 for every $\delta>0$ there exists
$\beta=\beta(D,\delta)>0$ such that 
if $H$ is generated according to the above model (with
an arbitrary choice of the $d(v)$'s), then
with probability $1-o(1)$ (as $m\to \infty$),
 every $S\subseteq A$ with $|S|\le\beta m$ has at least
$d(S)-(1+\delta)|S|$ neighbors in $B$.
\end{lemma}

\heading{Proof. } Let us write $w(S)=d(S)-(1+\delta)|S|$.
Since every $S$ has at least $|S|/D$ neighbors, it suffices
to consider only the $S$ with $w(S)\ge |S|/D$.
The calculation is a variation of that in \cite{HLW}.

For sets $S\subseteq A$, $|S|=s\le\beta m$, and
$W\subseteq B$, $|W|=\lfloor w(S) \rfloor$, let
$X_{S,W}$ be the event ``all neighbors of $S$ lie in $W$.''
If no $X_{S,W}$ occurs, then $H$ satisfies the conclusion
of the lemma. We have
$$
\Pr[X_{S,W}] = \frac{d(W)(d(W)-1)\cdots (d(W)-d(S)+1)}
{d(B)(d(B)-1)\cdots (d(B)-d(S)+1)}\le
\left(\frac{d(W)}{d(B)}\right)^{d(S)}\le
\left(\frac{D^2w(S)}{m}\right)^{d(S)}.
$$
Thus for $S$ fixed, the probability of $X_{S,W}$ occurring
for some $w$ is at most
$$
{|B| \choose w(S)} \left(\frac{D^2w(S)}{m}\right)^{d(S)}.
$$
Estimating the binomial coefficient as ${x\choose y}\le
(ex/y)^y$ and using $s/D\le w(S)\le Ds$, this can be bounded by
$$
\left(\frac{e|B|}{w(S)}\right)^{w(S)} 
\left(\frac{D^2w(S)}{m}\right)^{d(S)}
\le \left(\frac{e|B|}{s}\right)^{w(S)} 
\left(\frac{D^3 s}{m}\right)^{d(S)}\le 
C_1^s \left(\frac s m\right)^{(1+\delta)s}
$$
($C_1$ a constant independent of $\beta$).
Then the probability of any $X_{S,W}$ occurring at all is bounded
by 
$$
\sum_{1\le s\le \beta m} {|A| \choose s} C_1^s \left(\frac s m\right)^{(1+\delta)s}
\le \sum_s C_2^s \left(\frac s m\right)^{\delta s},
$$
where $C_2$ is another constant. The term for $s=1$
is $O(m^{-\delta})=o(1)$, and the ratio of consecutive
terms is at most $C_2\beta^\delta$, which can be
made smaller than $\frac 12$, say, by fixing
$\beta$ small enough.
Then the entire sum is $o(1)$ as claimed.
\proofend

\begin{corol} \label{c:almosttree}
In the setting of Lemma~\ref{sparse},
the following holds almost surely:
The average degree of the subgraph of $H$ induced by any
set of  at most $\beta m$ vertices is at most $2+2\delta$.
\end{corol}

\heading{Proof. } 
Let us consider the subgraph of $H$ induced by $S\cup T$,
$S\subseteq A$, $T\subseteq B$, $|S\cup T|\le \beta m$.
By Lemma~\ref{sparse}, we may
assume that $S$ has at least $w(S)$ neighbors.
Hence at least $w(S)-|T|$ neighbors of $S$ do not lie
in $T$, and each such neighbor ``consumes'' at least
one edge among the $d(S)$ edges incident to $S$.
Thus the number of edges in $H[S\cup T]$ is at
most $d(S)-w(S)+|T|= (1+\delta)|S|+|T|$, and the
average degree is at most $2+2\delta$
(actually, at most $2+\delta$, if we use symmetry
and assume $|S|\ge |T|$).
\proofend

\medskip

\heading{Proof of Theorem~\ref{t:6pluseps}. }
We let $\rho=\rho(\eps)>0$ be a sufficiently
small constant, and  we set the parameters 
of our random graph model as follows:
we let $d(v)=5$ for some $\rho m$ vertices in $A$,
and all remaining vertices in $A\cup B$ have
 $d(v)=4$. Clearly, the average degree of $H$
is $4+\Omega(\rho)$, and thus if we choose
$\delta$ sufficiently small in terms of $\rho$
and let $\beta=\beta(5,\delta)$, then 
Lemma~\ref{local_density}  and Corollary~\ref{c:almosttree}
guarantee $\mcc2(L(H))\ge \beta m$ almost surely.

The maximum degree of $L(H)$ is $7$; the degree-7 vertices
of the line graph
correspond to edges of $H$ incident to the
$\rho m$ vertices of degree $5$. It remains to bound
the maximum average degree of $L(H)$.

To this end, we apply Corollary~\ref{c:almosttree} once again,
this time with $\delta=\frac12$, say, and we let 
$\beta_0=\beta(5,\frac 12)$
be the corresponding parameter. We note that $\beta_0$
is independent of $\rho$, and hence we can assume 
that $\beta_0/\rho$ is sufficiently large.

Now let $F\subseteq E(H)$ be an arbitrary subset
of edges, and let $U$ be the set of all vertices incident
to edges of $F$.
If $|U|\le\beta_0 m$, then by Corollary~\ref{c:almosttree}
the graph $K:=(U,F)$ has average degree at most $2.5$. 
The average degree of $L(K)$ is 
$$
\overline d(L(K))=\frac 1{|F|} \sum_{\{u,v\}\in F}(\deg_{K}(u)+
\deg_{K}(v)-2)=
\frac1{|F|}\sum_{u\in U}\deg_{K}(u)(\deg_{K}(u)-1).
$$
If we denote by $x_i$ the fraction of vertices 
 $u\in U$ with $\deg_{K}(u)=i$, $i=1,2,\ldots,5$,
then the $x_i$ satisfy the constraints
 $\sum_{i=1}^5 x_i=1$ and
$\sum_{i=1}^5 ix_i\le 2.5$ (this reflects the 
bound on the average degree of $K$), and we have
$\overline d(L(K))
=2(2x_2+6x_3+12x_4+20x_5)/(x_1+2x_2+3x_3+4x_4+5x_5)$.
One can use, e.g.,  linear programming to verify that
the above constraints imply $\overline d(L(K))\le 6$
as needed.

If, on the other hand, $|U|>\beta_0 m$, then 
since $H$ has at most $\rho m$ vertices of degree $5$,
these vertices constitute at most
 $\rho/\beta_0$ fraction of $U$,
and all other vertices have degree at most $4$.
Hence at most a small fraction of the vertices
of $L(K)$ can have degree $7$, while others have
degrees at most $6$, and
it follows that the average degree of $L(K)$
can be pushed below $6+\eps$ by making $\rho$
sufficiently small. Theorem~\ref{t:6pluseps}
is proved.
\proofend

\section{The Hamming cube}

An interesting example, where
$\mcc2$ can be determined exactly, 
is the Hamming cube. Let $Q_d$ denote the $d$-dimensional Hamming cube
with vertex set $\{0,1\}^d$ and with two vectors
$u,v\in \{0,1\}^d$ adjacent in $Q_d$ if they differ
in exactly one coordinate. Let $L(Q_d)$ be the line
graph of $Q_d$.

\begin{prop}\label{p:cube} For every even $d$ we have
$\mcctwo(L(Q_d))=\frac d4 2^{d/2}$.
\end{prop}

\heading{Proof. } The upper bound is witnessed by the following
coloring: Color an edge $\{u,v\}\in E(Q_d)$ red if $u$ and $v$
differ in one of the first $d/2$ coordinates, and blue otherwise.
Then the monochromatic components are $(d/2)$-dimensional subcubes.

For the lower bound, it suffices to show that
whenever the edges of $Q_d$ are colored red and blue,
there exists a monochromatic connected subgraph
with at least $\frac d4 2^{d/2}$ edges.
Let us assume that, e.g., blue is the majority color;
that is, at least $\frac d4 2^{d}$ edges are blue.
Let $B_1,B_2,\ldots,B_k$ be the connected
components of the blue subgraph, and let $m_i$ be the number
of  blue edges in $B_i$. Hence $\sum_{i=1}^k m_i
\ge \frac 12|E(Q_d)|=\frac d4 2^{d}$.

We recall the following formulation of the
edge-isoperimetric inequality for the cube
(E.g., \cite{bol86}, Chapter~16):
 Every subgraph of $Q_d$
on $v$ vertices has at most $\frac 1 2v\log_2 v$ edges.
Let $\beta_i>0$ be the real number satisfying
$\frac 12\beta_i\log_2\beta_i=m_i$ (where $\beta_i=1$ for $m_i=0$).
Thus, $\beta_i$ is a lower bound for the number of vertices
of $B_i$, and consequently, $\sum_{i=1}^k \beta_i\le 2^d$.
Assuming for contradiction that $m_i<\frac d42^{d/2}$ for all $i$,
we have $\beta_i<2^{d/2}$ for all $i$, and thus
$$
\sum_{i=1}^k m_i =\sum_{i=1}^k \frac 12\beta_i\log_2\beta_i<
\frac d4\sum_{i=1}^k\beta_i\le \frac d4\sum_{i=1}^k |V(B_i)|
\le \frac d4 2^d.
$$
But as was noted above,
$\sum_{i=1}^k m_i\ge \frac d42^{d}$, and this contradiction
establishes the proposition.
\proofend

\medskip

The proof also shows that the monochromatic components
in any extremal coloring have to be $(d/2)$-dimensional subcubes.

\section{Open problems}

There are quite a few
interesting open questions suggested by the present paper.
Here are some of them.

\begin{enumerate}
\item How large can $\mcctwo(G)$ be for graphs of maximum degree $6$?
By the planar HEX lemma, the triangulated planar grid 
is an example with $\mcctwo(G)= \Theta( \sqrt n)$, but this
is at present the best we know.
\item A special case of the previous question,
which seems interesting in its own right, is when $G=L(H)$ 
for some $4$-regular $H$. The best lower bound
we know is $\Omega(\log n)$, from a construction
by Alon et al.\ \cite{ADO+03}, where $H$ is 
a $4$-regular graph of logarithmic girth.
\item The examples we know for planar graphs $G$ with large $\mcctwo(G)$
and $\mcc 3(G)$ have at least one vertex of high degree.
Can anything better be said if we assume that $G$ has bounded degrees?
More specifically, the following was asked in \cite{KMR+97}:
Is there a function $f$ such that for every planar graph $G$ of maximum degree
$\Delta$ we have $\mcc3(G)\le f(\Delta)$?
\item For two colors we cannot hope for constant monochromatic component
size in bounded-degree planar graphs, as shown
by a triangulated planar grid, but
similar to the previous question, we can
ask if there exists a function $g$ such
that every $n$-vertex planar graph $G$ of maximum degree $\Delta$
satisfies $\mcc2(G)\le g(\Delta)\sqrt n$.
\item There is still a gap between the best bounds we know for $\mcc t(G)$
for graphs from minor-closed families of graphs. Can this gap be closed?
\item A question suggested to us by Emo Welzl concerns the possible
behavior of $\mcc t(G)$ when the chromatic number of $G$ as well as its
number of vertices are known.
\item
The proof of Theorem~\ref{t:6pluseps}
can be adapted to show that for any fixed $t$ there
exist $n$-vetrex graphs with maximum degree $4t-1$ and with
 $\mcc t(G) = \Omega(n)$. Note that in \cite{ADO+03} it was shown
that $\mcc t(G)$ is not bounded by a constant 
even for graphs of maximum
degree $4t-2$; however, that proof gives only 
a logarithmic lower bound for
$\mcc t(G)$ (cf.\ our first open problem). From the other direction
\cite{HST03} show that $\mcc t(G)$ is bounded by a constant for
all $t$ and all graphs $G$ with maximum degree at most $3t-1$. 
The constant $3$ here is not optimal, since
the same paper shows that for some constant $\eps>0$ 
and all sufficiently large $t$, the
value $\mcc t(G)$ is bounded by a constant for
all graphs $G$ of maximum degree at most $(3+\eps)t$. It would be
interesting to find the asymptotic behaviour of the maximal value of $\mcc
t(G)$ for graphs $G$ with maximum degree $d$ in the intermediate range
$(3+\epsilon)t<d<4t-2$. In particular, it would be interesting to know if there
exist $t$ and $d$ for which the above maximum is sublinear but not a constant.
\item There are several natural conjectures
pertaining to $\mcc t$ for triangulations of the $d$-dimensional
grid graph. These questions suggest an interesting ``combinatorial
dimension theory'' waiting to be discovered. More on this subject
can be found in~\cite{MP}.
\end{enumerate}

\section*{Acknowledgements}

We would like to thank  Maria
Chudnovsky and Tibor Szab\'o for fruitful conversations and
to Eyal Ackerman and Robin Thomas for important references.

\end{document}